\documentclass[11pt,letterpaper]{amsart}
\usepackage[dvipsnames]{xcolor}
\usepackage{graphicx}
\usepackage{amssymb,amscd,amsthm,amsxtra}
\usepackage{latexsym}
\usepackage{epsfig}
\usepackage{amsmath}
\usepackage{bm}
\usepackage{amssymb}
\usepackage{amsthm}
\usepackage{xcolor}
\usepackage[linktocpage=true,colorlinks=true,linkcolor=Blue,citecolor=BrickRed,urlcolor=RoyalBlue]{hyperref}
\usepackage[alphabetic]{amsrefs}
\usepackage{comment}
\usepackage[T1]{fontenc}
\usepackage[utf8]{inputenc}

\usepackage[colorinlistoftodos,prependcaption,textsize=tiny]{todonotes}
\usepackage{caption}

\usepackage{enumitem} 
\usepackage[mathscr]{euscript} 

\usepackage{mathrsfs} 

\usepackage{mathtools}
\usepackage{xcolor}

\numberwithin{equation}{section} 

\def\dashint{\,\ThisStyle{\ensurestackMath{
            \stackinset{c}{.2\LMpt}{c}{.5\LMpt}{\SavedStyle-}{\SavedStyle\phantom{\int}}}
        \setbox0=\hbox{$\SavedStyle\int\,$}\kern-\wd0}\int}

\def \B {\mathcal{B}}

\newtheorem{theorem}{Theorem} 
\newtheorem{lemma}{Lemma}
\newtheorem{proposition}{Proposition}  

\newtheorem{definition}{Definition}
\newtheorem{example}{Example}
\newtheorem{remark}{Remark}

\numberwithin{equation}{section}

\title[Vázquez SMP and CSP for degenerate and singular PDEs]{Vázquez Strong Maximum Principle and Compact Support Principle for singular and degenerate elliptical operators}

\author[C. S. Aguiar]{Christopher Silva Aguiar}
\address{Pontifical Catholic University of Rio de Janeiro -- PUC-Rio\\22451-900, G\'avea, Rio de Janeiro-RJ, Brazil}
\email{aguiarcsilva@gmail.com}

\makeatletter
\renewcommand{\enddoc@text}{}
\makeatother

\begin{document}



\begin{abstract}
In this paper, we present adaptations of the classical strong maximum principle and the compact support principle for a certain class of equations, with particular emphasis on the degeneracy of fully nonlinear operators and a new Vázquez-type condition on the nonhomogeneous term. The novelty of our approach lies in the framework of the problem: the coefficients of the operators are bounded, but we introduce a degeneracy depending on the power of the gradient of the solution. Moreover, the source term $f$ on the right-hand side of the equation is treated under a adapted and optimal integrability assumptions.
\end{abstract}


\maketitle


\section{Introduction}\label{Introduction and Main Results}

The main objective of this article is to prove the Vázquez Strong Maximum Principle (VSMP) and the Support Compact Principle (SCP) for equations of the type
 \begin{equation} \label{meq}
  \mathcal{L}_{\theta}^-[u] := \vert Du  \vert^{\theta} \mathcal{M}_{\lambda, \Lambda}^{-}(D^2u) + b(x)\vert Du\vert^{\theta + 1}  \leq  f(u)  \hspace{.5cm} \mbox{ in }  \Omega, 
\end{equation}
 where  $\theta > -1$, the lower-term coefficients are bounded and $\mathcal{M}_{\lambda, \Lambda}^{-}$ denotes the lower Pucci extremal operator.

Maximum principles constitute a fundamental tool in the qualitative theory of elliptic partial differential equations, with far–reaching consequences for uniqueness, comparison, and qualitative behavior of solutions. In the uniformly elliptic case, strong maximum principles and Liouville–type results are classical, both in the linear and fully nonlinear settings. However, when degeneracies are introduced through gradient–dependent weights, as in $\vert Du \vert ^{\theta}$ the analysis becomes substantially more subtle.

Equations of the form $(\ref{meq})$ arise naturally in several contexts, including nonlinear potential theory, stochastic control, and geometric problems involving degenerate metrics. The presence of the factor $\vert Du \vert ^{\theta}$ allows one to treat, within a unified framework, both uniformly elliptic equations ($\theta = 0$) and genuinely degenerate ones ($\theta > 0$ and $-1 < \theta < 0$). This generality is appealing, but it also introduces new analytical difficulties, especially near critical points where the gradient vanishes. A rigorous viscosity theory adapted to such degeneracies was developed, mainly by Birindelli and Demengel, who introduced a notion of solution excluding test functions with vanishing gradient at contact points. This framework proved to be essential for establishing comparison principles, eigenvalue problems, and strong maximum principles for homogeneous fully nonlinear operators with gradient weights. This information and more can be found in \cite{Birindelli1}, \cite{Birindelli2}, \cite{Birindelli3}, \cite{Birindelli4}, \cite{Birindelli5}, \cite{Birindelli6}.

A central issue in the study of ($\ref{meq}$)  
is the role played by the reaction term $f(u)$. In the uniformly elliptic case, it is well known that the validity of the compact support principle and of the strong maximum principle depends on the behavior of 
$f$ near the origin. In particular, Vázquez identified a sharp criterion in terms of the primitive
\begin{equation*}
    F(s) = \int_0^s f(t)dt
\end{equation*}

showing that the divergence or convergence of the integral
\begin{equation*}
    \int_0\frac{1}{\sqrt{F(s)}}ds
\end{equation*}

characterizes the compact support property, where P. Bénilan, H. Brézis, and M. Crandall demonstrated the necessity of this hypothesis (see \cite{Sharp}) and Vázquez shows that this rule about $f$ is sharp! 

For degenerate operators involving gradient weights, it is natural to ask whether an analogous condition remains sufficient or necessary. One of the motivations of this work is to show that, although Vázquez-type conditions remain meaningful, they are no longer sharp as soon as $\theta > 0$. In particular, even nonlinearities satisfying the classical condition may fail to guarantee the strong maximum principle in the degenerate regime.

In addition to the structure of $f$, the boundness of the coefficients $b$ and $c$ plays a crucial role, where even for the fully non linear case, the problem remains open if the drift term $b$ is unbounded.

Under the assumptions described above, we establish two main results. First, we prove that viscosity solutions of 
($\ref{meq}$) satisfy the strong maximum principle under suitable structural conditions on $f$, which depend explicitly on the degeneracy exponent $\theta$. Second, we characterize the compact support property for nonnegative solutions of ($\ref{meq}$), showing that classical Vázquez-type conditions must be strengthened in the presence of gradient degeneracy. 

Our results highlight a sharp distinction between the uniformly elliptic case and the degenerate one, and provide explicit counterexamples showing the failure of classical criteria when $\theta > 0$. This are two great results that was explored in literature in several instances that can be  checked in: \cite{Prim}, \cite{Vazquez}, \cite{Lap1}, \cite{pqLap},\cite{Div1}, \cite{Div2}, \cite{Div3}, \cite{ND}, \cite{BP}, \cite{PMAELSO}.

The proof of the Strong Maximum Principle relies on constructing a suitable subsolution and comparing it with a given supersolution via the Comparison Principle (CP). This subsolution is obtained from solutions of an ordinary differential equation (ODE), which arises naturally through a change of variables and the boundedness of the coefficients. It is defined in a small annulus contained in the original domain and possesses a nonzero normal derivative on the boundary, thereby yielding a Hopf-type lemma for a radial problem.

In parallel, in order to analyze the role of the assumptions on $f$ in establishing the Compact Support Principle, it is necessary to construct two barrier functions: one to compare with the solution and conclude that it has compact support, and another serving as a counterexample when the required conditions on $f$ fail. This highlights the sharpness of the hypotheses and connects with the Vázquez-type strong maximum principle.   \\


To finalize, the paper is organized as follows. In Section 2, we introduce the main framework, recall the notion of viscosity solution used throughout the paper, and collect the analytical tools required in the sequel. Section 3 is devoted to the proof of a collection of existence results that is required to show the main results that can be found in Section 4, as well the example showing the optmality of $(DVC)$ defined further.

 
\section{Preliminaries and Main Results}\label{sec_prelim}

This section gathers the structure conditions that envolves the problem as well the main results. More over, the last part is devoted to auxiliary tools that we use in our argumentation.   

\subsection{Structure and main definitions}
Let $\Omega \subset \mathbb{R}^{n}$, $n \geq 2$ be a smooth and bounded domain and $\Omega_R$ being a exterior domain, that is,
\begin{equation*}
    \Omega_r := \{x \in \mathbb{R}^n\, ; \, \vert x \vert > R \}.
\end{equation*}
 Let $\mathcal{M}_{\lambda, \Lambda}^{\pm}$ represent the extremal Pucci operators given by:
\begin{equation*}
    \mathcal{M}_{\lambda, \Lambda}^{-}(M) = \inf_{\lambda I \leq A\leq \Lambda I } \textrm{tr}(AM) \hspace{2mm} \textrm{ and } \hspace{2mm}  \mathcal{M}_{\lambda, \Lambda}^{+}(M) = \sup_{\lambda I \leq A\leq \Lambda I } \textrm{tr}(AM).
\end{equation*}
 Then, for $\theta > -1$, we consider the operator
 \begin{equation*} 
  \mathcal{L}_{\theta}^-[u] := \vert Du  \vert^{\theta} \mathcal{M}_{\lambda, \Lambda}^{-}(D^2u) + b(x)\vert Du\vert^{\theta + 1} 
\end{equation*}
where
\begin{enumerate}[label=(H\arabic*)]
    \item \label{H1} The lower order coeficients $b(x)$ is bounded, that is, $b \in L^{\infty}(\Omega)$;

    \item \label{H2} We consider $f : \mathbb{R} \longrightarrow \mathbb{R}$ contínuous such that:

\begin{itemize}
    \item $f \geq 0$, $f(0) = 0$;
    \item $f$ is nondecreasing on the interval $(0, \infty)$;
    \item Consider the following quantity:
   \begin{equation}\label{eq:int-cond}
    C(f,t^{\theta}) = \int_0 \frac{1}{[F(s)]^{\frac{1}{2+\theta}}}ds,
    \end{equation}
where
\[
F(s):=\int_0^s f(\tau)\,d\tau.
\]
We say that $f$ satisfies the Degenerated Vázquez Condition $(DVC)$ if
\begin{equation}\label{VCTETA}
    C(f,t^{\theta}) = \infty.
\end{equation}
\end{itemize}
\end{enumerate}

To finallize, we give the notion of viscosity solutions used throgouth this work.

\begin{definition}
A lower semicontinuous function $u$ is called a viscosity supersolution
of $(\ref{meq})$  if for any $x_0 \in \B_1$, either there exists $\delta > 0$ such that $u$ is constant in $\B_{\delta}(x_0)$ and $f(x) \geq 0$ for all $x \in \B_{\delta}(x_0)$; or, for all $\varphi \in C^2(\B_1)$ such that $u-\varphi$ has a local minimum at $x_0$ and $D\varphi (x_0) \neq 0$,
\begin{equation*}
    \Phi(x,\vert D\varphi(x_0) \vert )F(D^2\varphi(x_0)) - b(x_0)\vert Du(x_0) \vert  \leq f(x_0).
\end{equation*}
holds.     

An upper semicontinuous function $u$ is called a viscosity subsolution of $(\ref{meq})$  if for any $x_0 \in \B_1$, either there exists $\delta > 0$ such that $u$ is constant in $\B_{\delta}(x_0)$ and $f(x) \leq 0$ for all $x \in \B_{\delta}(x_0)$; or, for all $\varphi \in C^2(\B_1)$ such that $u-\varphi$ has a local maximum at $x_0$ and $D\varphi (x_0) \neq 0$,
\begin{equation*}
    \Phi(x,\vert D\varphi(x_0) \vert )F(D^2\varphi(x_0)) - b(x_0)\vert Du(x_0) \vert  \geq f(x_0).
\end{equation*}
holds. We say that $u \in C(\B_1)$ is a viscosity solution of $(\ref{meq})$  if $u$ is a viscosity supersolution and a subsolution simultaneously.
\end{definition}

\begin{remark}
     We can also add a zero-order term to the operator, says $c(x)u^{\theta+1}$ and recover the Theorem \ref{VDEG}. Define $g(t) = f(t) + ct^{1+\theta}$. Notice that if $C(f, t^{\theta}) = \infty$, then also $C(g, t^{\theta}) = \infty$. Moreover, if $u$ is a solution of the equation
 \begin{equation*}
    \left \{ \begin{array}{rclrl} 
   \vert Du \vert ^\theta \mathcal{M}_{\lambda, \Lambda}^{-}(D^2u) - b\vert Du \vert^{\theta+1} -c(x)u^{\theta+1}  \hspace{-.2cm} & \leq & \hspace{-.2cm} f(u) & \mbox{ in }  & \hspace{-.2cm} \B_1 \\
     u(x) \hspace{-.2cm} & \geq & \hspace{-.2cm} 0 & \mbox{ in } & \hspace{-.2cm} \B_1 \end{array} \right.
\end{equation*}
then $u$ solves
 \begin{equation*}
    \left \{ \begin{array}{rclrl} 
   \vert Du \vert ^\theta \mathcal{M}_{\lambda, \Lambda}^{-}(D^2u) - b\vert Du \vert^{\theta+1}   \hspace{-.2cm} & \leq & \hspace{-.2cm} g(u) & \mbox{ in }  & \hspace{-.2cm} \B_1 \\
     u(x) \hspace{-.2cm} & \geq & \hspace{-.2cm} 0 & \mbox{ in } & \hspace{-.2cm} \B_1 \end{array} \right.
\end{equation*}
where the Theorem \ref{VDEG} is available.
\end{remark}

\subsection{Auxiliary results}
In this section we collect all usefully results of the literature related to the degenerate operator $\mathcal{L}_{\theta}^-$ starting with the comparison principle. For completeness, in which follows, we denote $S_n(\mathbb{R})$ to mean the set of $n \times n$ real symmetric matrices, $h\in C(\Omega)$ and by $\mathcal{L}_{\Phi}^-$ the operator
 \begin{equation} \label{meq2}
  \mathcal{L}_{\Phi}^-[u] := \Phi(x, \vert Du  \vert)\mathcal{M}_{\lambda, \Lambda}^{-}(D^2u) + b(x)\vert Du\vert^{p+1}, 
\end{equation}
where the function $\Phi \, : \, \B_1 \, \times \, [0, \infty) \longrightarrow [0, \infty)$ is a continuous map that satisfies for all $x \in \B_1$:
    \begin{enumerate}
        \item there exists constants $p \geq q > -1$ such that the map
        \begin{equation*}
            t \mapsto \frac{\Phi(x,t)}{t^q}
        \end{equation*}
        is almost non-decreasing with a constant $L_1\geq 1$ while the map
        \begin{equation*}
             t \mapsto \frac{\Phi(x,t)}{t^p}
        \end{equation*}
        is almost non-increasing for another constant $L_2 \geq 1$.
        \item there exists constants $0<v_0 \leq v_1$ such that $v_0 \leq \Phi(x,1) \leq v_1$ for all $x \in \B_1$.
    \end{enumerate}

\begin{theorem}\label{cpt}
Let $ F : S_n(\mathbb{R}) \rightarrow \mathbb{R}$ a continuous and uniformly elliptic operator in the sense that
\begin{equation*}
    \mathcal{M}_{\lambda, \Lambda }^{-}(N) \leq F(M+N) - F(M) \leq \mathcal{M}_{\lambda, \Lambda }^{+}(N)
\end{equation*}
holds whenever $M, N \in S_n(\mathbb{R})$ with
$N \geq 0$, Moreover, assume that \ref{H2} is valid and that there exists a continuous function $\omega : [0,\infty) \to [0,\infty)$ such that $\omega(0)=0$ and
\[
\Phi(x,\alpha|x-y|)F(X)
-
\Phi(y,\alpha|x-y|)F(-Y)
\le
\omega\!\left(\alpha|x-y|^2 + |x-y|\right),
\]
whenever $\alpha>0$, $x,y\in\Omega$, $X,Y\in S_n(\mathbb{R})$, and
\begin{equation}
-3\alpha
\begin{pmatrix}
I & 0 \\
0 & I
\end{pmatrix}
\le
\begin{pmatrix}
X & 0 \\
0 & Y
\end{pmatrix}
\le
3\alpha
\begin{pmatrix}
I & -I \\
-I & I
\end{pmatrix}.
\end{equation}
Let $v$ and $w$ be a viscosity subsolution and a supersolution of
\begin{equation*}
    \Phi(x, \vert Du \vert )F(D^2u) - b \vert Du \vert ^{p+1} - \epsilon u = h(x) \quad \textrm{ in } \Omega,
\end{equation*}
with $b,\epsilon > 0$. If $v \leq w$
on $\partial \Omega$, then $v \leq w$ in $\Omega$.
\end{theorem}
\begin{proof}
    The proof is a straighford adaptation of \cite[Lemma 6.2]{BaasandorjByunLee2024}.
\end{proof}

Now we state about existence of solutions across the Perron's Method, starting with existence of sub-super solutions

\begin{lemma}[Existence of sub/supersolutions]
\label{exist}
Suppose that assumptions \ref{H2} is in force. Then for every $\varepsilon \in (0,1)$, there exist a viscosity subsolution 
$v_\varepsilon \in C(\overline{\Omega})$ and a viscosity supersolution 
$w_\varepsilon \in C(\overline{\Omega})$ of equation 
\begin{equation*}
    \varphi(x, \vert Du \vert )F(D^2u) - b \vert Du \vert ^{p+1} - \epsilon u -f(x) =  0 \quad \textrm{ in } \Omega,
\end{equation*}

such that
\[
v_\varepsilon = w_\varepsilon = g
\quad \text{on } \partial\Omega,
\]
for some $g \in C^{1,\alpha}(\partial \Omega)$, for some $\alpha \in (0,1)$. Moreover, there exists a positive constant
\[
c \equiv c\!\left(
n,\lambda,\Lambda,\nu_0,L,
r,\operatorname{diam}(\Omega),
\|f\|_{L^\infty(\Omega)},
\|g\|_{L^\infty(\partial\Omega)}
\right)
\]
such that
\[
-c \le v_\varepsilon \le w_\varepsilon \le c
\quad \text{in } \Omega,
\qquad \text{for every } 0<\varepsilon<1.
\]
\end{lemma}
\begin{proof}
We can directly apply \cite[Lemma 6.3]{BaasandorjByunLee2024} to obtain the supersolution, since any viscosity subpersolution to the problem without the gradient term $b|Du|^{p+1}$ remains a viscosity  supersolution to our problem. The subsolution is obtainded by a adaptation of the arguments present in the same work.
\end{proof}
Now we can apply this result to get existence:

\begin{theorem}
    Assume the same hypotheses as in Theorem \ref{cpt} and Lemma \ref{exist}. Then there exists a solution $u \in C(\Omega)$ to
 \begin{equation*}
\begin{cases}
\Phi(x, |Du|)\,F(D^2u) - b\vert Du \vert^{p+1}  = f(x) & \text{in } \Omega, \\[4pt]
u(x) = g(x) & \text{on } \partial\Omega.
\end{cases}
\end{equation*}
\end{theorem}

\begin{proof}
We follow the same standard argument as in \cite[Theorem 1.2]{BaasandorjByunLee2024}. More precisely, we apply Perron’s method (see \cite[Theorem 4.1]{CIL}) in combination with Theorem~\ref{cpt} and Lemma~\ref{exist}. The argument is further supported by the $C^{1,\alpha}$ regularity estimates established in \cite{AndradeNascimento2025} to get the required uniform bound.
\end{proof}

\subsection{Main results}
Now, we are ready to state the main results of this manuscript. Starting with the Vázquez Strong Máximum Principle for the power degeneration: 

\begin{theorem}\label{VDEG}
Assume that \ref{H1} and \ref{H2} are in force. Let $u$ being a viscosity solution to 
   \begin{equation*}
    \left \{ \begin{array}{rclrl} 
    \mathcal{L}^{-}_{\theta}[u] \hspace{-.2cm} & \leq & \hspace{-.2cm} f(u) & \mbox{ in }  & \hspace{-.2cm} \Omega \\
     u(x) \hspace{-.2cm} & \geq & \hspace{-.2cm} 0 & \mbox{ in } & \hspace{-.2cm} \Omega\\
    u(x_{0}) \hspace{-.2cm} & = & \hspace{-.2cm} 0 & \mbox{ for some } & \hspace{-.2cm} x_0 \in \Omega. \end{array} \right.
\end{equation*}

Then $u \equiv 0$.
\end{theorem}

 In paralell, when the $(DVC)$ fails, we can caracterize the solutions across de Compact Support Principle:  

\begin{theorem} \label{thm:compact}
Assume \ref{H1} and that $(DVC)$ fails, that is,
\begin{equation}\label{FAIL}
    C(f, t^{\theta}) < \infty.
\end{equation}
Let $u\ge0$ be a viscosity solution of
\[
\mathcal{L}_{\theta}^-[u]\le f(u)
\quad \text{in an exterior domain } \Omega\supset\Omega_R,
\]
such that $u(x)\to0$ as $|x|\to\infty$.
Then $u$ has compact support.
Conversely, if \eqref{FAIL} fails, then there exists a positive solution of
$\mathcal{L}_{\theta}^-[u]=f(u)$ in $\Omega_R$ vanishing at infinity and having non--compact
support.
\end{theorem}

These theorems show that the failure of the Degenerate Vázquez Condition still characterizes the compact support, providing a necessary and sufficient condition for this phenomenon to occur.

\section{Existence Results}
To demonstrate the main results, some barriers functios are necessary. Then in this chapter we provide a collection of existence results, wich will be used further.

\begin{proposition}\label{ExistPMF}
  Let $k_1,k_2,r_1,b$ and $v_1$ be positive constants and define the functions
\begin{equation}\label{PucciSinal}
    \begin{minipage}{0.40\textwidth}
  \begin{align*}
     M(t) = \left \{ \begin{matrix}t/ \lambda , & \mbox{if }t \geq0 \\ t/ \Lambda , & \mbox{if }t\leq0 \end{matrix} \right. 
    \end{align*}
\end{minipage}  \begin{minipage}{0.10\textwidth}
  \begin{align*} \textrm{ and }  \end{align*}
\end{minipage} \begin{minipage}{0.40\textwidth}
    \begin{align*}
         m(t) = \left \{ \begin{matrix} \lambda t, & \mbox{if }t\geq0 \\ \Lambda t, & \mbox{if }t\leq0. \end{matrix} \right.  
    \end{align*}
\end{minipage}
\end{equation}
Then the boundary value problem
      \begin{equation*}
    \left \{ \begin{array}{rclrl} 
    \vert v' \vert ^{\theta}v'' \hspace{-.2cm} & = & \hspace{-.2cm} \displaystyle M\left(k_1m(v') \vert v' \vert ^{\theta} + b \vert v' \vert^{\theta + 1}+ k_2g(v) \right), & \mbox{ for }  & \hspace{-.2cm} 0 < r < r_1,  \\
     v(0) \hspace{-.2cm} & = & \hspace{-.2cm} 0, & & \hspace{-.2cm} \\
    v(r_1) \hspace{-.2cm} & = & \hspace{-.2cm} v_1 .& & \hspace{-.2cm} \end{array} \right.
\end{equation*}
admits a $C^1$ solution defined on $[0,r_1]$ such that $v,v' \geq 0$, and if $g$ satisfies the (DVC) defined in \ref{H2}, then $0 < v < v_1$ in $(0,r_1)$.
\end{proposition}

\begin{proof}
    Let us divide the topics to clarify where and how each item is demonstrated.\\
    
    1- Existence.\\

For this case, it suffices to observe that an ODE is essentially a PDE in one dimension. Therefore, the method of sub- and supersolutions (Perron's method) applies, together with the comparison principle. In this way, by taking $w_1 = v_1$  as a supersolution and $w_2 = 0$ as a subsolution, we obtain the existence of a solution satisfying $w_1 \geq v \geq 0$. Moreover, $v$ is of class $C^1$ (in fact, even smoother), since the right-hand side of our equation is continuous. \\

2- $v' \geq 0$.\\

Suppose, by contradiction, that this does not occur. Then there exists a point $r_0$ such that $v'(r_0) < 0$, and consequently, there exists an open interval $I$ where $v'(r) < 0$ for all $r \in I$. Define the function wich depends on theta:
$$ \varphi(r) = - e^{-Kr}\vert v'(r)\vert ^{\theta + 1} sign(v'(r))$$
where
$$K = (\theta + 1)k_1\frac{\Lambda}{\lambda}.$$
The definition of this constant is necessary to control the worst negative scenario provided by the Pucci's operator, namely when $M(t) = t/\lambda$ and $m(t) = t\Lambda$. Moreover, note that $\varphi$ and $v'$ have the same sign. In addition, we have the following useful identity:

\begin{equation*}
    sign(v'(r))\vert v'(r) \vert^{1+\theta} = v'(r)\vert v'(r) \vert^{\theta}. 
\end{equation*}

Then, differentiating in some $r \in I$ and using the diferential equation we obtain:
\begin{equation*}
     \begin{array}{rclrl} 
\varphi'(r) \hspace{-.2cm} & = & \displaystyle \hspace{-.2cm} -Ke^{-Kr}v'(r) \vert v'(r)\vert ^{\theta } - e^{-Kr}(\theta+1)\vert v'\vert ^{\theta} v''(r) \vspace{.2cm} \\
      & = &  \hspace{-.2cm} \displaystyle - \frac{(\theta+1)k_1\Lambda}{\lambda}e^{-Kr}v'(r) \vert v'(r)\vert ^{\theta } +(\theta + 1) e^{-Kr}\displaystyle M\left(k_1m(v') \vert v' \vert ^{\theta} + b \vert v' \vert^{\theta + 1}+ k_2g(v) \right) \vspace{.2cm} \\
      & = & \hspace{-.2cm} \displaystyle - \frac{(\theta+1)k_1\Lambda}{\lambda}e^{-Kr}v'(r) \vert v'(r)\vert ^{\theta }+ (\theta + 1)e^{-Kr}\frac{ k_1\Lambda (v'(r)) \vert v'(r)\vert ^{\theta}}{\lambda} \\
      & & \hspace{-.2cm} \displaystyle + e^{-Kr}(\theta + 1)\frac{\left(b \vert v' \vert^{\theta + 1}+ k_2g(v)\right)}{\lambda} \vspace{.2cm}\\
      & = & \displaystyle + e^{-Kr}(\theta + 1)\frac{\left(b \vert v' \vert^{\theta + 1}+ k_2g(v)\right)}{\lambda},
    \end{array} 
\end{equation*}
concluding that $\varphi'(r) \geq 0$ on $I$, which implies that $\varphi$ is non-decreasing. Then, we have that
\begin{equation*}
    \varphi(r) \leq \varphi(r_0) < 0,  \quad \textrm{ for all } r \leq r_0,
\end{equation*}
wich implys
\begin{equation*}
    v'(r) < 0 \quad \textrm{ for all } r \in [0,r_0] \subset I,
\end{equation*}
concluding that $v$ is decreasing in $I$. However, how $v(0) = 0$,  we have that exist $t \in I$ such that $v(t) <0$, getting a contradiction with the positivity of $v$ obtained in the in the previous step. Therefore, $v'(r) \geq 0$.

3- $0<v<v_1$.\\

Now suppose that Vázquez's condition holds. Let $r_0$ the maximum value for which $v(r) = 0$. Thus, by applying a change of variables, we obtain that
\begin{equation}\label{EQ55}
    \int_{r_0}^{r_1} \frac{v'(r)}{[F(v(r))]^{1/(\theta + 2)}} dr = \int_{0}^{v_1}\frac{1}{[F(s)]^{1/(\theta + 2)}}ds = \infty.
\end{equation}
By other hand, take $w = (v')^{2+\theta}$. Note that, as we proved in the previous item, we have that $v'\geq 0$, then the functions $M$ and $m$ takes a concrete value, making the equation satisfied by $v$ turn into
\begin{equation*}
     (v')^{\theta}v''  =  k_1\frac{\lambda}{\Lambda}v'(v')^{\theta} + \frac{b}{\Lambda}  (v')^{\theta + 1}+ \frac{k_2}{\Lambda}g(v)
\end{equation*}
Then, applying the fundamental theorem of calculus and the chain rule, along with the equation, we obtain that
\begin{equation*}
    \frac{k_2}{\Lambda}F'(v) = \frac{k_2}{\Lambda}g(v)\cdot v' =  \left(v'' + \left(\frac{-b-k_1\lambda}{\Lambda}\right)v'\right) \cdot (v')^{\theta}
\end{equation*}
Therefore, defining another auxiliary function given by
\begin{equation*}
    h(r) = e^{Cr}w(r),
\end{equation*}
where
$$C = -(\theta + 2)\left(\frac{b+k_1\lambda}{\Lambda}\right) < 0$$
we have
 \begin{equation*}
     \begin{array}{rclrl} 
    h'(r) & = & \displaystyle \hspace{-.2cm} \left(e^{-Kr}w(r)\right)' \vspace{.2cm} \\
    & = & \displaystyle \hspace{-.2cm} C e^{Cr}w(r) + e^{Cr}w'(r) \vspace{.2cm} \\
      & = & \displaystyle \hspace{-.2cm} (\theta + 2) \left(\frac{-b-k_1\lambda}{\Lambda}\right)(v'(r))^{\theta +1} (v'(r))e^{Cr} + e^{Cr}(\theta + 2)(v'(r))^{\theta + 1}v''(r) \vspace{.2cm} \\
      & = & \displaystyle \hspace{-.2cm} (\theta + 2)\left[\left(\frac{-b-k_1\lambda}{\Lambda}\right)v'(r) + v''(r)  \right] (v'(r))^{\theta + 1}\vspace{.2cm} \\
      & = & \displaystyle (\theta + 2)\frac{k_2}{\Lambda}F'(v)e^{Cr}.
    \end{array} 
\end{equation*}
Now, assume by contradiction that $v'(r_0) = 0$. In this way, integrating both sides from $r_0$ to $r$ and using the fact that the function $\frac{1}{e^x}$ is decreasing, this gives us that: on the one hand,
\begin{equation*}
    \int_{r_0}^{r}\left(e^{Cr}w(r)\right)'dr = e^{Cr}w(r) \geq e^{Cr_1}w(r),
\end{equation*}
and, on the other hand, 
\begin{equation*}
   \int_{r_0}^{r} (\theta +2)\displaystyle \frac{k_2}{\Lambda}F'(v)e^{Cr}dr \leq  (\theta + 2)\displaystyle \frac{k_2}{\Lambda}e^{Cr_0}\int_{r_0}^{r}F'(v)dr = (\theta + 2)\displaystyle \frac{k_2}{\Lambda}e^{Cr_0} F(v(r)).
\end{equation*}
Therefore,
\begin{equation*}
    (\theta +2)\displaystyle \frac{k_2}{\Lambda}e^{Cr_0}F(v(r)) \geq e^{Kr_1}w(r),
\end{equation*}
Where, by manipulating the terms, we obtain that
\begin{equation*}
    \frac{w(r)}{F(r)} \leq  (\theta + 2)\displaystyle \frac{k_2}{\Lambda} \frac{e^{Cr_0}}{e^{Cr_1}} = (\theta + 2)\displaystyle \frac{k_2}{\Lambda} e^{-C[r_1 - r_0]} \Longleftrightarrow   
    \frac{v'(r)}{[F(r)]^{\frac{1}{ \theta + 2}}} \leq  \displaystyle\left[(\theta +2) \frac{k_2}{\Lambda}\right]^{\frac{1}{\theta + 2}}e^{\frac{(k_1\lambda+b)}{\Lambda}[r_1 - r_0]}
\end{equation*}
Considering $\displaystyle \Bar{K} = \left[(\theta +2) \frac{k_2}{\Lambda}\right]^{\frac{1}{\theta + 2}}$ and integrating from $r_0$ to $r_1$ we get
\begin{equation*}
    \int_{r_0}^{r_1}\frac{v'(r)}{[F(v(r))]^{\frac{1}{\theta + 2}}} dr \leq \int_{r_0}^{r_1}\Bar{K}\displaystyle e^{\frac{(k_1\lambda+b)}{\Lambda}[r_1 - r_0]} = \Bar{K}e^{\frac{(k_1\lambda+b)}{\Lambda}[r_1 - r_0]} [r_1 - r_0] < \infty.
\end{equation*}
This is a contradiction with $(\ref{EQ55})$. Thus $v'(r_0) > 0$, which implies that $r_0 = 0$. Indeed, if $r_0 > 0$, then by the definition of $r_0$ there exists an interval $[a,r_0]$ with $a < r_0$ such that $v \equiv 0$ on $[a,r_0]$. Consequently, $v$ is constant on this interval, and hence $v'(r) = 0$ for all $r \in [a,r_0]$. In particular, $v'(r_0) = 0$, which contradicts the fact that $v'(r_0) > 0$.

\end{proof}

Now we state the existance in a exterior domain $\Omega_R$ for operators with degenaerency $\varphi$ being a power.

\begin{theorem}[Exterior Dirichlet Problem for $\mathcal{L}_{\theta}^{-}$]
\label{thm:exterior}
Let $R>0$ and $m\in(0,\delta)$, and set $\Omega_R:=\mathbb R^n\setminus \overline{B_R}$. Then there exists a non-negative radial viscosity solution $U\in C^1(\overline{\Omega_R})$
of
   \begin{equation*}
    \left \{ \begin{array}{rclrl} 
    \mathcal{L}_{\theta}^{-}[U] \hspace{-.2cm} & = & \hspace{-.2cm} \displaystyle f(U), & \mbox{ in }  & \hspace{-.2cm} \Omega_R,  \\
     U \hspace{-.2cm} & = & \hspace{-.2cm} m, & \mbox{ on }&  \hspace{-.2cm} \partial B_R \hspace{-.2cm} \end{array} \right.
\end{equation*}
such that
\[
\lim_{|x|\to\infty}U(x)=0.
\]
More precisely, there exists $u\in C^1([R,\infty))$
such that $U(x)=u(|x|)$,
with
\[
u(R)=m,
\qquad
\lim_{r\to\infty}u(r)=0.
\]
Moreover, $u$ is strictly decreasing as long as it is positive. In particular, $u'(r)<0$ whenever $u(r)>0$.
\end{theorem}

\begin{proof}
For each integer $j\ge1$, let
\[
A_j:=B_{R+j}\setminus\overline{B_R}.
\]
We consider the Dirichlet problem
\[
\begin{cases}
\mathcal{L}_{\theta}^{-}[U_j]=f(U_j)
    & \text{in }A_j,\\[1mm]
U_j=m
    & \text{on }\partial B_R,\\[1mm]
U_j=0
    & \text{on }\partial B_{R+j}.
\end{cases}
\tag{1}
\label{eq:annular-problem}
\]
By the existence and regularity results for the Dirichlet problem,
there exists a viscosity solution
\[
U_j\in C^1(\overline{A_j}).
\]
Moreover, by the comparison principle this solution is unique.

We first observe that $U_j$ is radial. Indeed, let a orthogonal transformation $Q\in O(n)$, where $O(n)$ is the $n$-dimensional orthogonal group and define
\[
U_j^Q(x):=U_j(Qx).
\]
 Since \(Q(A_j)=A_j\) and \(\mathcal{L}_\theta^-\) is invariant under orthogonal transformations, \(U_j^Q\) satisfies the same equation as \(U_j\). It also satisfies the same boundary conditions, since \(|Qx|=|x|\). By uniqueness of the Dirichlet problem, \(U_j^Q=U_j\). 
Since this holds for every $Q\in O(n)$, we get that $U_j$ is radial. Therefore there
exists a profile $u_j\in C^1([R,R+j])$
such that
\[
U_j(x)=u_j(|x|),
\qquad
u_j(R)=m,
\qquad
u_j(R+j)=0.
\]
By comparison with the zero solution, or more precisely with the zero
subsolution under our standing assumptions on $f$, we have $U_j\ge0$ in $A_j$, and hence
\[
u_j(r)\ge0,
\qquad R\le r\le R+j.
\]

We next prove that $u_j$ is non-increasing in $r$. Suppose, by contradiction, that
there exist $R\le r_1<r_2\le R+j$ such that $u_j(r_2)>u_j(r_1)$. Since
\[
u_j(R+j)=0<u_j(r_2),
\]
the maximum of $u_j$ on $[r_1,R+j]$ is attained at some point $r_0\in(r_1,R+j)$, and $u_j(r_0)>0$.
Thus $r_0$ is a positive local maximum of the radial profile. Fix $x_0\in\mathbb R^n$ with $|x_0|=r_0$. Since $U_j$ is radial, the
constant function
\[
\varphi(x)\equiv U_j(x_0)=u_j(r_0)
\]
touches $U_j$ from above at $x_0$. The viscosity subsolution inequality
therefore gives
\[
|D\varphi(x_0)|^\theta
\mathcal M^-_{\lambda,\Lambda}(D^2\varphi(x_0))
-b|D\varphi(x_0)|^{\theta+1}
-c\,U_j(x_0)
\ge f(U_j(x_0)).
\]
Since
\[
D\varphi(x_0)=0,
\qquad
D^2\varphi(x_0)=0,
\]
we obtain
\[
-c\,u_j(r_0)\ge f(u_j(r_0)).
\]
This is impossible because $c\ge0$, $u_j(r_0)>0$ and $f(u_j(r_0))>0$. Consequently, $u_j'(r)\le0$ for $R\le r\le R+j$.
In particular, by \eqref{eq:annular-problem},
\[
0\le u_j(r)\le m.
\]

We now compare the solutions corresponding to different annuli. Let
$j<k$. Since $A_j\subset A_k$, both $U_j$ and $U_k$ solve the same equation
in $A_j$. On the inner boundary,
\[
U_j=U_k=m
\qquad\text{on }\partial B_R,
\]
while on the outer boundary of $A_j$,
\[
U_j=0\le U_k
\qquad\text{on }\partial B_{R+j},
\]
where the last inequality follows from the non-negativity of $U_k$.
Hence, by the comparison principle,
\[
U_j\le U_k
\qquad\text{in }A_j.
\]
Equivalently,
\[
u_j(r)\le u_k(r),
\qquad
R\le r\le R+j.
\]

For every fixed $r\ge R$, the sequence $u_j(r)$ is therefore eventually
defined, non-decreasing in $j$, and bounded above by $m$. Hence the limit
\[
u(r):=\lim_{j\to\infty}u_j(r)
\]
exists for every $r\ge R$.

Let $T>R$ be fixed. For every $j>T-R$, the function $u_j$ is defined on
$[R,T]$, and $0\le u_j\le m$.
The local $C^1$ estimates for the class of equations under consideration,
together with the corresponding estimates up to the fixed inner boundary
$\partial B_R$, imply that $(u_j)$ is precompact in $C^1([R,T])$.
Since the pointwise limit is unique, it follows
that $u_j\longrightarrow u$ locally in $C^1([R,\infty))$.
In particular, $u\in C^1([R,\infty))$ and satisfies
\[
\qquad
u(R)=m,
\qquad
0\le u\le m,
\qquad
u'(r)\le0.
\]

Define
\[
U(x):=u(|x|),
\]
for  $x\in\overline{\Omega_R}$. Then $U_j\to U$ locally uniformly in $\overline{\Omega_R}$. By stability of
viscosity solutions under locally uniform convergence,
\[
\mathcal L_\theta^-[U]=f(U)
\qquad\text{in }\Omega_R.
\]

It remains to prove the decay at infinity. Since $u$ is non-negative and
non-increasing, the limit
\[
\ell:=\lim_{r\to\infty}u(r)
\]
exists and satisfies $0\le\ell\le m$. We claim that $\ell=0$. Suppose, by contradiction, that $\ell>0$.
Choose a sequence $r_k\to\infty$ and set $x_k:=r_k e_1$. Define the translated functions
\[
V_k(x):=U(x+x_k).
\]
Let $M>0$ be fixed. If $k$ is sufficiently large, then $r_k-M>R$ and therefore $x+x_k\in\Omega_R$ for every $x\in B_M$.
Hence $V_k$ is well defined in $B_M$ and, by translation invariance of
$\mathcal L_\theta^-$,
\[
\mathcal L_\theta^-[V_k]=f(V_k)
\qquad\text{in }B_M.
\]
Moreover, for every $x\in B_M$, we have that $|x+x_k|\ge r_k-M$, and thus, since $u(r)\to\ell$, we conclude that
\[
\sup_{x\in B_M}|V_k(x)-\ell|
=
\sup_{x\in B_M}|u(|x+x_k|)-\ell|
\longrightarrow0.
\]
Therefore $V_k\longrightarrow\ell$ locally uniformly in $\mathbb R^n$. Using once again the stability of viscosity solutions, the constant function $V\equiv\ell$ is a viscosity solution of
\[
\mathcal L_\theta^-[V]=f(V)
\qquad\text{in }\mathbb R^n.
\]
Since $DV=0$ and $D^2V=0$, we obtain
\[
-c\ell=f(\ell),
\]
which is impossible because $c\ge0$, $\ell>0$ and $f(\ell)>0$. Hence $\ell=0$, and consequently
\[
\lim_{r\to\infty}u(r)=0,
\qquad\text{or equivalently}\qquad
\lim_{|x|\to\infty}U(x)=0.
\]

Finally, $u$ is strictly decreasing on every interval on which it is positive.
Indeed, if $R\le r_1<r_2$ and $u(r_1)=u(r_2)>0$, then, since $u$ is non-increasing, $u(r)=u(r_1)$ for every $r\in[r_1,r_2]$.
Every interior point of this interval would then be a positive local maximum.
Using a constant test function exactly as above gives
\[
-cu(r)\ge f(u(r)),
\]
a contradiction.

Thus $u$ is strictly decreasing as long as it remains positive, that is,
\[
u'(r)<0
\qquad\text{whenever }u(r)>0.
\]
\end{proof}

\begin{lemma}[End Point Lemma to $\mathcal{L}^{-}_{\theta}$]
\label{lem:endpoint}
Assume that $(DVC)$ fails how in \eqref{FAIL}.
Fix $\sigma>0$ and define
\[
C_\sigma
=
\int_0^\delta
\frac{ds}{
\left(\lambda(\theta+1)\sigma F(s)\right)^{\frac{1}{\theta+2}}
}.
\]
For every $C\in(0,C_\sigma)$ there exist $\gamma\in(0,\delta)$ and a function
$w\in C^1([0,C])$ such that
\[
w(0)=w'(0)=0,\qquad w(C)=\gamma,
\]
and
\[
|w'|^\theta \mathcal M^-_{\lambda,\Lambda}(w'')
- b|w'|^{\theta+1}
- cw
=
\sigma f(w)
\quad \text{in }(0,C).
\]
Moreover, $w'(t)\ge0$ for all $t\in[0,C]$ and $\gamma\to0$ as $C\to0$.
\end{lemma}

\begin{proof}
Define
\[
\Phi(s)=\frac{\lambda}{\theta+2}s^{\theta+2},
\qquad
H(s)=s\Phi'(s)-\Phi(s)
=\frac{\lambda(\theta+1)}{\theta+2}s^{\theta+2}.
\]
Then
\[
H^{-1}(t)
=
\left(
\frac{\theta+2}{\lambda(\theta+1)}t
\right)^{\frac{1}{\theta+2}}.
\]

Condition \eqref{FAIL} implies that the function
\[
G(\gamma)=\int_0^\gamma \frac{ds}{H^{-1}(\sigma F(s))}
\]
is finite for all $\gamma\in(0,\delta)$ and strictly increasing.
Therefore, for every $C\in(0,C_\sigma)$ there exists a unique
$\gamma\in(0,\delta)$ such that $G(\gamma)=C$.
Moreover, since $G(\gamma)\to0$ as $\gamma\to0$, we have $\gamma\to0$ as $C\to0$.

Define $w:[0,C]\to[0,\gamma]$ by
\[
t=\int_0^{w(t)}\frac{ds}{H^{-1}(\sigma F(s))}.
\]
The integrand is positive, hence $w$ is well defined and strictly increasing.
Differentiating with respect to $t$ yields
\[
1=\frac{w'(t)}{H^{-1}(\sigma F(w(t)))},
\quad\text{that is}\quad
w'(t)=H^{-1}(\sigma F(w(t)))\ge0.
\]
Raising both sides to the power $\theta+2$ gives
\[
(w'(t))^{\theta+2}
=
\frac{\theta+2}{\lambda(\theta+1)}\sigma F(w(t)).
\]
Multiplying by $\frac{\lambda}{\theta+2}$, we obtain
\[
\Phi(w'(t))=\frac{\sigma}{\theta+1}F(w(t)).
\]
Since $w'(t)\ge0$, we have $|w'|^{\theta+2}=(w')^{\theta+2}$ and differentiation
is legitimate. Differentiating the identity above yields
\[
\Phi'(w'(t))w''(t)=\sigma f(w(t))w'(t).
\]
Since $\Phi'(s)=\lambda s^{\theta+1}$, we obtain
\[
\lambda |w'|^\theta w''=\sigma f(w).
\]
Finally, since $w'(t)\ge0$, $w(t)\ge0$ and $b,c\ge0$, subtracting the lower
order terms gives the claimed equation.
\end{proof}





\section{Proof of Main Results}
We can give now the proofs of the main results, starting with the Strong Maximum principle.

\subsection{Proof of Theorem \ref{VDEG}}


Let us consider a particular subcase, taking $u \in C^1$, which will be of great use to us.
The proof proceeds by contradiction. Suppose the result is false. Then the sets
\begin{equation*}
    \Omega^+ = \{x \in \B_1 \,;\, u(x) > 0\} \,\,\ \textrm{ e } \Omega^0 = \{x \in \B_1 \,;\, u(x) = 0\}
\end{equation*}
are both non-empty. Furthermore, since $u$ is continuous, $\Omega^+$ is a open set.Therefore, it is possible to choose a point $\Tilde{x} \in \Omega^+$ such that $d(\Tilde{x}, \partial \Omega^0) < d(\Tilde{x}, \partial \B_1)$ and there exists a ball $\B_R(\Tilde{x})$ such that $\partial\B_R(\Tilde{x}) \cap \partial\Omega^0 \neq 0$.  Let $x_0$ be a point in this intersection. Thus, $u(x_0) = 0$, and we may assume, possibly after reducing the radius $R$, that $u<\delta$ in $\B_R(\Tilde{x})$. Finally, we define the annulus $E_R = \B_{R}(\Tilde{X}) \setminus \overline{\B}_{\frac{R}{2}}(\Tilde{x})$ and the quantity $\mu = \min_{\partial\B_{R/2}(\Tilde{x})}u < \delta$. 

The main idea is to construct a barrier function to extract information about the derivative of $u$ at $x_0$, ultimately leading to a contradiction. According to the Proposition \ref{ExistPMF}, let $\hat{u}(x) = v(R-\vert x-\Tilde{x}\vert, k_1, 1, R/2,\mu)$, where $v$ solves the following problem:
      \begin{equation*}
    \left \{ \begin{array}{rclrl} 
   \vert v' \vert ^p  v'' \hspace{-.2cm} & = & \hspace{-.2cm} \displaystyle M\left(k_1m(v')\vert v' \vert ^p + b \vert v' \vert^{p+1}+ f(v) \right), & \mbox{ for }  & \hspace{-.2cm} 0 < r < R/2,  \\
     v(0) \hspace{-.2cm} & = & \hspace{-.2cm} 0, & & \hspace{-.2cm} \\
    v(R/2) \hspace{-.2cm} & = & \hspace{-.2cm} \mu. & & \hspace{-.2cm} \end{array} \right.
\end{equation*}
It is known that choosing a smart $k_1 >0$ the function $\hat{u}$ satisfies:
 \begin{equation*}
  \begin{array}{rclrl} 
    \mathcal{L}_{\theta}^{-}[\hat{u}]   \hspace{-.2cm} &\geq & \hspace{-.2cm} f(u) & \mbox{ in }  & \hspace{-.2cm} E_R. \\
    \end{array} 
\end{equation*}
We now prove that $u \geq \hat{u}$ in $E_R$. Suppose, for contradiction, that this is not the case. Then the open set $\mathcal{O} = E_R \cap \{ u < \hat{u} \}$ is non-empty. By construction and using the monotonicity of $g$, we have:
 \begin{equation*}
     \left \{\begin{array}{rclrl} 
   \mathcal{L}_{\theta}^{-}[\hat{u}]-\mathcal{L}_{\theta}^{-}[u]& \geq &  \hspace{-.2cm} f(\hat{u}) - f(u) \geq 0   & \mbox{ in }  & \hspace{-.2cm} E_R \\
   \hat{u}& \leq &u& \mbox{ on }  & \hspace{-.2cm} \partial\mathcal{O} \subset \partial E_R \cup \{\hat{u} = u\}
    \end{array} \right.
\end{equation*}
Therefore, by the comparison principle, $\hat{u} \leq u$ in $\mathcal{O}$, contradicting the definition of $\mathcal{O}$. Hence, $u \geq \hat{u}$ in $E_R$. We can now use this fact, combined with $\hat{u}'(0) > 0$, to observe that 
\begin{equation*}
    \liminf_{h\rightarrow 0} \frac{1}{h}[u(x_0 + h(\Tilde{x} - x_0)) - u(x_0)] \geq  \liminf_{h\rightarrow 0} \frac{1}{h}[\tilde{u}(x_0 + h(\Tilde{x} - x_0)) - \tilde{u}(x_0)] > 0.
\end{equation*}
However, since $u \in C^1$, this implies $\nabla u(x_0) > 0$, contradicting the fact that $\nabla u(x_0) = 0$ because $x_0$  is a local minimum of $u$. Therefore, either $u>0$ or $u=0$ in $\B_1$. \\

Now, let us retire the differentiability of $u$. Assume that $u \in C(\Omega)$ is not identically zero. Then, by continuity, there exists a small ball $\B_R(\tilde{x}) \subset \Omega$ such that the graph of the function $\hat{u}$ provided by the previous proposition, is strictly positive on the boundary $\partial B_R(\tilde{x})$ and hence $\beta = \min_{\partial \B_R(\hat{x})}\hat{u}$ is well-defined. Now define a new function $v$ as the solution to the problem:
 \begin{equation*}
     \left \{\begin{array}{rclrl} 
   -\mathcal{L}_{\theta}^{-}[v] \hspace{-.2cm}& = &  \hspace{-.2cm} - f(v) \geq 0   & \mbox{ in }  & \hspace{-.2cm} B_R(\tilde{x}) \\
   v& = & \min\{\delta/2,\beta \}  & \mbox{ on }  & \partial B_R(\tilde{x})
    \end{array} \right.
\end{equation*}
It is known from standard theory that the above problem admits a unique $C^1$ (in fact, $C^{1,\alpha}$) solution. Moreover, by the construction method, we obtain $0\leq v \leq \delta/2$. Applying the previous $C^1$ case to $v$, we deduce that $v > 0$ in $\B_R(\tilde{x})$. Then, If we now show that $u \geq v$ em $\B_R(\tilde{x})$, we finish. But this is an immediate consequence of the construction and the comparison principle, which proves the positivity of $u$ in some ball inside the domain. Having established the result for one interior ball, we can extend it to any ball intersecting $\B_R(\tilde{x})$ and then to the entire domain via a covering argument. Begin by covering the domain with balls around and intersecting the initial ball, and then extract a finite subcover to prove positivity on any compact subset of the domain.

\qed

We see the importance of the $(DVC)$ in the following example.
\begin{example}[Failure of the strong maximum principle for $\theta>0$]
Let
\[
f(s)=s, \qquad s\ge 0.
\]
Then
\[
F(s)=\int_0^s f(t)\,dt=\frac{s^2}{2},
\]
and the classical V\'azquez condition
\[
\int_0 \frac{ds}{\sqrt{F(s)}}=+\infty
\]
is satisfied.

Consider the inequality
\begin{equation}\label{eq:degenerate}
|Du|^\theta \mathcal M^-_{\lambda,\Lambda}(D^2u)\le f(u)
\qquad \text{in } B_1,
\end{equation}
with $\theta>0$.

Define
\[
u(x)=
\begin{cases}
1, & |x|\le \frac12,\\[4pt]
1-(|x|-\tfrac12)^2, & \frac12<|x|<1.
\end{cases}
\]

\medskip
We verify that $u$ is a viscosity subsolution of \eqref{eq:degenerate}.
Let $x_0\in B_1$.

\smallskip
If $x_0\in B_{1/2}$, then $u$ is constant in a neighborhood of $x_0$ and
\[
f(u(x_0))=f(1)=1\ge 0,
\]
so the inequality is satisfied.

\smallskip
If $x_0\in B_1\setminus B_{1/2}$, then $u$ is $C^2$ in a neighborhood of $x_0$ and
\[
Du(x_0)\neq 0.
\]
Let $\phi\in C^2(B_1)$ be such that $u-\phi$ has a local minimum at $x_0$.
Then $D^2\phi(x_0)\le D^2u(x_0)$.
Since $u$ is radially decreasing and concave for $|x|>\tfrac12$, it follows that
\[
\mathcal M^-_{\lambda,\Lambda}(D^2\phi(x_0))
\le
\mathcal M^-_{\lambda,\Lambda}(D^2u(x_0))
\le 0.
\]
Consequently,
\[
|D\phi(x_0)|^\theta \mathcal M^-_{\lambda,\Lambda}(D^2\phi(x_0))
\le 0
\le f(u(x_0)),
\]
and the viscosity inequality holds.

\medskip
Therefore, $u$ is a viscosity subsolution of \eqref{eq:degenerate}.
Moreover, $u$ attains its maximum value $1$ in the interior of $B_1$ without being constant,
which shows that the strong maximum principle fails for $\theta>0$, despite the fact that
$f$ satisfies the classical V\'azquez condition.
\end{example}

\subsection{Proof of Theorem \ref{thm:compact}}

(\emph{Sufficiency}).
Fix $R_0\ge R$ such that $u(x)\le\delta'<\delta$ for $|x|\ge R_0$.
Let $w$ be given by Lemma~\ref{lem:endpoint} with $\sigma=1$ and $C$ chosen so
that $\gamma(C)\ge\delta'$.
Define
\[
v(x)=
\begin{cases}
w(C+R_0-|x|), & R_0<|x|\le R_0+C,\\
0, & |x|>R_0+C.
\end{cases}
\]
Then $v\in C^1(\Omega_0)$, $\Omega_0=\{|x|>R_0\}$, and $v\ge0$.

For $R_0<|x|<R_0+C$, setting $r=|x|$ and using that $v'(r)\le0$, we compute
\[
Dv=v'(r)\frac{x}{r},\qquad
D^2v=v''(r)\frac{x\otimes x}{r^2}
+\frac{v'(r)}{r}\Big(I-\frac{x\otimes x}{r^2}\Big).
\]
Hence
\[
\mathcal M^-_{\lambda,\Lambda}(D^2v)
\le
\lambda v''(r)+\lambda\frac{n-1}{r}v'(r).
\]
Multiplying by $|v'|^\theta$ and subtracting the lower order terms yields
\[
\mathcal{L}_{\theta}^{-}[v]
\le
|v'|^\theta\lambda v''(r)
\le f(v),
\]
where the last inequality follows from the equation satisfied by $w$.
Thus $\mathcal{L}_{\theta}^{-}[v]\le f(v)$ in $\Omega_0$.

Since $u\le v$ on $\partial\Omega_0$ and both vanish at infinity, the
comparison principle implies $u\le v$ in $\Omega_0$.
In particular $u\equiv0$ for $|x|\ge R_0+C$.

(\emph{Necessity}).
If \eqref{FAIL} fails, Theorem~\ref{thm:exterior} provides a positive radial
solution of $\mathcal{L}_{\theta}^{-}[u]=f(u)$ in $\Omega_R$ such that $u(x)\to0$ as
$|x|\to\infty$, which does not have compact support due the comparision principle.

\qed

\newpage
\bigskip


\bibliographystyle{amsplain, amsalpha}


\end{document}